\documentclass[reqno, 12pt]{amsart}

\usepackage[utf8]{inputenc}
\usepackage{amsmath, amsthm, amssymb, amsfonts}
\usepackage{mathrsfs, latexsym, amscd}
\usepackage{graphicx, xcolor}
\usepackage{hyperref}
\usepackage{booktabs, makecell, multirow}

\usepackage[margin=1in]{geometry}
\numberwithin{equation}{section}

\title[Exact Scan Statistics for $k = 3$ and $k = N-1$]{Exact Closed-Form Formulae for Linear and Circular Continuous Scan Statistics: $P_c(N - 1; N, w)$, $P_c(3; N, w)$, and $P(3; N, w)$}
\author{Haowei Yuan}
\address{Peking University}
\email{hwyuan25@stu.pku.edu.cn}

\subjclass[2020]{62E15, 62G30, 60D05, 05A15}
\keywords{Scan statistics, Circular spacings, Geometric probability, Exact distribution, Generating functions}

\begin{document}

\begin{abstract}
  The continuous linear $P(k; N, w)$ and circular scan statistics $P_c(k; N, w)$ are fundamental tools in probability and spatial statistics, frequently used to detect clustering in uniform data. Let $X_1, X_2, \dots, X_N$ be independently and uniformly distributed random variables on a unit interval or unit ring. The exact distribution of these scan statistics relies on the minimum window width required to capture exactly $k$ points. Furthermore, the survival function $1 - P_c(k; N, w)$ directly corresponds to the geometric probability that if $N$ arcs of length $1 - w$ are uniformly and randomly placed on a unit circle, every point on the circle is covered at least $N + 1 - k$ times. Historically, evaluating the exact cumulative distribution functions, $P(k; N, w)$ and $P_c(k; N, w)$, relies heavily on complex recursive approximations. In this paper, we bypass these traditional recursive methods to derive direct, generalized closed-form expressions for some linear and circular continuous scan statistics. Specifically, we present the exact analytical solutions for $P_c(N - 1; N, w)$, $P_c(3; N, w)$, and $P(3; N, w)$ for arbitrary values of $N$ and window width $w$. These newly derived closed-form expressions not only provide exact baseline distributions for extreme spacings but also significantly simplify computational complexity compared to existing iterative approaches.
\end{abstract}

\maketitle

\section{Introduction}

Let $X_1, X_2, \dots, X_N$ be a sequence of $N$ random variables independently and uniformly distributed on the unit interval, with their order statistics denoted as $X_{(1)} \le X_{(2)} \le \dots \le X_{(N)}$. The linear continuous scan statistic relies on finding the minimum window width required to capture exactly $k$ points, defined algebraically by the $k$-step spacing $W(k) = \min_{1 \le i \le N - k + 1} (X_{(i + k - 1)} - X_{(i)})$. To account for periodic boundaries in circular geometries, the corresponding minimum circular window width on a unit ring is defined as $W_c(k) = \min\left(W(k), \min_{N - k + 2 \le i \le N}(X_{(i + k - 1 - N)} + 1 - X_{(i)})\right)$. Consequently, the exact distributions of the linear and circular scan statistics can be formulated as the cumulative distribution functions of these spacing extremes: $P(k; N, w) = \mathbb{P}(W(k) \le w)$ and $P_c(k; N, w) = \mathbb{P}(W_c(k) \le w)$. This notation follows \cite{glaz2001scanningexact}.

Geometrically, the survival function $1 - P_c(k; N, w)$ is equivalent to the classic multiple coverage problem of a circle, often referred to as Stevens' problem \cite{stevens1939solution}. In this geometrical context, it represents the probability that $N$ uniformly placed random arcs of length $1 - w$ cover all points on the circumference at least $N + 1 - k$ times \cite{glaz1979multiple}.

Deriving exact probability distributions for continuous scan statistics is historically notorious for its mathematical difficulty. As a sliding window moves continuously across the domain, the event counts in adjacent, overlapping windows are highly correlated. Furthermore, the exact cumulative distribution function is not a smooth, unified curve; rather, it is a complex piecewise sequence of polynomials whose active domains change depending on the magnitude of $w$.

The calculation of exact distributions of $P(k; N, w)$ and $P_c(k; N, w)$ has a long history. \cite{burnside1928theory} found the exact closed-form expression for $P(N; N, w)$; \cite{parzen1960modern} found the exact closed-form expression for $P(2; N, w)$; \cite{naus1963clustering} found the exact closed-form expression for $P(N - 1; N, w)$. Note that the statistic $P(N - 1; N, w)$ differs significantly from $P_c(N - 1; N, w)$, because the linear case only has two constraints. \cite{naus1963clustering, naus1965distribution} derived a simple formula for $P(k; N, w)$ for $k > (N + 1)/2$. Due to the complexity of overlapping dependencies, determinantal approaches \cite{karlin1959coincidence, huntington1975simpler} struggled to produce manageable algebraic formulae for larger $N$ and arbitrary $w$, shifting much of the theoretical focus toward asymptotic approximations \cite{glaz2001scanningbounds}. Notably, \cite{vanElteren1961} calculated the exact closed-form expression for $P(3; N, w)$ with $N \le 8$ via direct integration.

In recent decades, the demand for exact calculations has largely been met through recursive algorithms and Markov chain imbedding techniques based on spacing dualities \cite{huffer1997computing, huffer1997approximating, huffer2001computing}, and an approach specially designed for circular cases \cite{lin2017two}. While these algorithmic approaches are exact, they often lack the analytical clarity of closed-form expressions and can suffer from exponential computational time complexities depending on the parameterization of $k$ and $N$. To date, unified, explicit algebraic representations for arbitrary variables have remained elusive.

In this paper, we bypass algorithmic recursion and determinantal matrices entirely to present a unified, purely algebraic framework for deriving the exact, closed-form piecewise polynomials of continuous scan statistics for all valid arbitrary parameters of $N$ and $w$.

The remainder of this paper is organized as follows. Section 2 establishes the fundamental sequence constraints, linking the continuous scan statistic parameters to spatial Lebesgue measures and constructs the recursive convolution framework. Section 3 transitions into the Laplace transform domain, and resolves the foundational generating function via a Riccati equation. Section 4 employs the Principle of Inclusion and Exclusion to adapt these functions to circular rings and linear boundary conditions. Finally, Section 5 applies the Lagrange Inversion Theorem to extract the exact binomial coefficients, then Section 6 presents the fully simplified, explicit closed-form formulae for $P_c(N-1; N, w)$, $P_c(3; N, w)$ and $P(3; N, w)$.

\section{Geometric Modeling and Recursive Formulations}

Let $f_n(x)$ denote the Lebesgue measure of the set
$\{(x_1, x_2, \dots, x_n) \in \mathbb{R}^n : 0 \le x_i \le 2, \sum_{i=1}^n (x_i - 1) = x, x_i + x_{i + 1} \le 2 \text{ for } 1 \le i \le n-1\}$.

Let $A$ be the set of indices $i$ such that $x_i \le 1$, and let $B$ be the set of the remaining indices (where $1 < x_i \le 2$).
The condition $x_i + x_{i + 1} \le 2$ implies two key properties: first, there are no adjacent indices in $B$; second, if an index $i \in B$ is adjacent to an index $j \in A$, then $|x_i - 1| \le |x_j - 1|$.
Consequently, if $n \ge 2$, the maximum value of $|x_i - 1|$ across all $i$ must be achieved by an index in $A$.

Next, when $n \ge 2$, we derive a recursive formula for $f_n(x)$ by conditioning on the index that maximizes $|x_i - 1|$. Let us enumerate the index $i$ at which this maximum occurs. As established, this index must belong to set $A$, meaning $x_i \le 1$. We define the maximum absolute difference as $p = 1 - x_i$, and we integrate over all possible values of $p \in [0, 1]$.

By fixing the maximum at index $i$, we partition the remaining sequence into two independent segments of lengths $i - 1$ and $n - i$. We then apply a geometric scaling transformation by dividing all shifted variables $(x_j - 1)$ by $p$, so the constraint $|x_j - 1| < p$ for every index $j \ne i$ becomes $0 \le x_j \le 2$, which aligns with the definition of $f(x)$. Prior to scaling, the sum of the remaining $n - 1$ shifted variables is $x - (x_i - 1) = x + p$. After dividing by $p$, the sum constraint for the scaled variables becomes $\frac{x + p}{p} = \frac{x}{p} + 1$.

Because the original system contains $n$ variables constrained by both a fixed sum and a fixed maximum, the subspace of the remaining variables has $n - 2$ dimensions. According to the properties of Lebesgue measure, dividing every dimension of this $(n-2)$-dimensional space by $p$ scales the overall measure by a factor of $p^{2 - n}$, so we must multiply by $p^{n - 2}$ to recover the original measure. Since the segments to the left and right of index $i$ are independent aside from their total sum, the measure of their combination is exactly the convolution of their individual measures. Thus, for $n \ge 2$, we obtain the recursive formula:

\begin{equation}
  f_n(x) = \sum_{i = 1}^n \int_0^1 (f_{i - 1} \ast f_{n - i}) \left(\frac{x}{p} + 1\right) \cdot p^{n - 2} \, dp.
\end{equation}

To complete the base cases for $n < 2$, we define $f_0(x)$ as the identity element of the convolution operator (the Dirac delta function), and $f_1(x)$ as the rectangular function corresponding to the allowed domain bounds:

\begin{equation}
  f_0(x) = \delta(x), \quad f_1(x) = \operatorname{rect}\left(\frac{x}{2}\right).
\end{equation}

\section{Bilateral Laplace Transform and Riccati Equation}

We use the bilateral Laplace transform to simplify the recursive convolution equation. Although there exist $x < 0$ for which $f_n(x) \neq 0$, the negative domain of $x$ is strictly bounded. Thus, we can treat it analogously to the ordinary Laplace transform. Defining the transform as $\mathcal{L}\{f\}(s) = \int_{-\infty}^{\infty} e^{-sx} f(x) \, dx$, we introduce a transformation $\tilde{f}_n(s) = s^n \cdot \mathcal{L}\{f_n(x)\}(s)$.

Evaluating the integral transformation, the recursive formula becomes:

\begin{equation}
  \tilde{f}_n(s) = \sum_{i = 1}^n \int_0^s e^p \cdot \tilde{f}_{i - 1}(p) \tilde{f}_{n - i}(p) \, dp.
\end{equation}

Differentiating with respect to $s$ yields:

\begin{equation}
  \frac{d}{ds} \tilde{f}_n(s) = e^s \sum_{i = 1}^n \tilde{f}_{i - 1}(s) \tilde{f}_{n - i}(s).
\end{equation}

Based on our established initial conditions, the base cases in the transform domain evaluate to:

\begin{equation}
  \tilde{f}_0(s) = 1, \quad \tilde{f}_1(s) = e^s - e^{-s}.
\end{equation}

To resolve this system for all $n$, we define the formal generating function $F(t, s) = \sum_{n=0}^{\infty} \tilde{f}_n(s) \cdot t^n$. Differentiating $F(t, s)$ with respect to $s$ and accounting for the initial condition $\tilde{f}_1(s)$, we obtain a Riccati differential equation:

\begin{equation}
  \frac{d}{ds} F(t, s) = t e^s \cdot F^2(t, s) + t e^{-s}.
\end{equation}

We solve this Riccati equation subject to the boundary condition $F(t, 0) = 1$. The roots of the associated characteristic equation are $\alpha_{1, 2} = \frac{1 \pm \sqrt{1 - 4t^2}}{2}$. Using the standard substitution $F(t, s) = -\frac{e^{-s}}{t} \frac{d}{ds} \ln Q(t, s) = -\frac{e^{-s}}{t} \frac{P(t, s)}{Q(t, s)}$, the solution takes the form:

\begin{equation}
  Q(t, s) = (\alpha_2 + t) e^{\alpha_1 s} - (\alpha_1 + t) e^{\alpha_2 s},
\end{equation}

\begin{equation}
  P(t, s) = \frac{d}{ds} Q(t, s) = \alpha_1 (\alpha_2 + t) e^{\alpha_1 s} - \alpha_2 (\alpha_1 + t) e^{\alpha_2 s}.
\end{equation}

To further simplify these coefficients, we recognize the structural connection to the Catalan generating function.
Let $z = \sqrt{1 - 4t^2}$ and define $C = \frac{1 - z}{2t^2}$. Consequently, we can express our parameters purely in terms of $C$ and $t$:

\begin{equation}
  z = \frac{2}{C} - 1, \quad \alpha_1 = \frac{1}{C} = 1 - t^2 C, \quad \alpha_2 = t^2 C = 1 - \frac{1}{C}.
\end{equation}

This yields the highly compact identities:

\begin{equation}
  \frac{\alpha_2}{\alpha_1} = (tC)^2, \quad \frac{\alpha_2 + t}{\alpha_1 + t} = tC.
\end{equation}

\section{Circular Chains and Principle of Inclusion-Exclusion}

Next, we consider the circular chain case corresponding directly to $P_c(N-1; N, w)$. Let $a_n(x)$ be the Lebesgue measure of the set satisfying $x_i \geq 0$, $\sum_{i=1}^n (x_i - 1) = x$, $x_i + x_{i + 1} \leq 2$, and the cyclic constraint $x_n + x_1 \leq 2$.

By uniquely breaking the circular symmetry at the index $i$ that maximizes $|x_i - 1|$ (which must belong to $A$), breaking the cycle leaves a linear sequence of length $n-1$:

\begin{equation}
  a_n(x) = n \cdot \int_0^1 f_{n-1}\left(\frac{x}{p} + 1\right) \cdot p^{n - 2} \, dp.
\end{equation}

Applying the scaled bilateral Laplace transform $\tilde{a}_n(s) = s^n \cdot \mathcal{L}\{a_n(x)\}(s)$ yields $\frac{\tilde{a}_n(s)}{n} = \int_0^s e^p \cdot \tilde{f}_{n - 1}(p) \, dp$. Defining the generating function $A(t, s) = \sum_{n=2}^{\infty} \frac{\tilde{a}_n(s)}{n} \cdot t^n$, we substitute and integrate to establish a link to the linear function $F(t, p)$:

\begin{equation}
  A(t, s) = \int_0^s t e^p \cdot F(t, p) \, dp = \int_0^s -\frac{d}{dp} \ln Q(t, p) \, dp = -\ln\left(\frac{Q(t, s)}{Q(t, 0)}\right),
\end{equation}

\begin{equation}
  \tilde{a}_n(s) = -n [t^n] \ln\left(\frac{Q(t, s)}{Q(t, 0)}\right).
\end{equation}

Evaluating the denominator yields $Q(t, 0) = -z$. Thus, $\tilde{a}_n(s) = n \cdot [t^n] \left\{ -\ln\left(\frac{Q(t, s)}{-z}\right) \right\}$.

For the specific cluster thresholds $P_c(3; N, w)$ and $P(3; N, w)$, we apply the Principle of Inclusion and Exclusion (PIE). Let $b_n(x)$ and $c_n(x)$ be the respective Lebesgue measures modeling the cyclic and linear chains where the constraint is $x_i + x_{i + 1} \geq 2$. By applying PIE, we partition this into components where adjacent sums are strictly less than 2. Removing a constraint severs the connection, meaning the total measure is the convolution of partitioned chains, with each remaining connection contributing a coefficient of $-1$. Note that in this case, when $n = 1$, the value of the only variable is unconstrained.

To adjust $F(t, s)$ for the $n=1$ unconstrained case, we define a function $R(t, s)$:

\begin{equation}
  \frac{R(t, s)}{Q(t, s)} = F(t, s) + t e^{-s} \implies R(t, s) = \alpha_2 (\alpha_2 + t) e^{-\alpha_1 s} - \alpha_1 (\alpha_1 + t) e^{-\alpha_2 s}.
\end{equation}

Using standard combinatorial relationships of generating functions for weighted rings and chains, we evaluate the transforms:

\begin{equation}
  \tilde{b}_n(s) = n (-1)^n [t^n] \left\{ A(t, s) - \ln\left(1 + \left(\frac{R(t, s)}{Q(t, s)} - 1\right)\right) \right\} = n (-1)^n [t^n] \left\{ -\ln\left(\frac{R(t, s)}{-z}\right) \right\},
\end{equation}

\begin{equation}
  \tilde{c}_n(s) = [t^{n+1}] (-1)^{n-1} \cdot \frac{1}{1 + \left(\frac{R(t, s)}{Q(t, s)} - 1\right)} \cdot (t e^s)^2 = (-1)^{n-1} [t^{n-1}] \left\{ \frac{Q(t, s)}{R(t, s)} \cdot e^{2s} \right\}.
\end{equation}

\section{Inverse Laplace Transform and Coefficient Extraction}

Applying the shifted inverse bilateral Laplace transform yields the identity, where $H(x)$ is the Heaviside step function:

\begin{equation}
  \mathcal{L}^{-1}\{e^{ns} s^{-m}\}(x) = \frac{(x + n)^{m - 1}}{(m - 1)!} \operatorname{H}(x + n).
\end{equation}

For generating functions containing terms of the form $e^{c z s}$ (where $c$ is a constant and $z$ is a formal power series in $t$ with $z(0) = 1$), we formally expand the operator using a Taylor series approach:

\begin{equation}
  \mathcal{L}^{-1}\{e^{czs} \cdot s^{-p}\}(x) = \mathcal{L}^{-1}\left\{\sum_{i = 0}^{\infty} \frac{(c(z - 1))^i s^i}{i!} e^{c s} s^{-p}\right\}(x) = \frac{(x + c z)^{p-1}}{(p-1)!} \operatorname{H}(x + c).
\end{equation}

Similarly, for the logarithmic terms essential to $a_n$ and $b_n$, the expansion yields:

\begin{equation}
  \mathcal{L}^{-1}\{-\ln(1 - a e^{czs}) \cdot s^{-p}\}(x) = \sum_{i = 1}^{\infty} \frac{a^i}{i} \frac{(x + i c z)^{p - 1}}{(p - 1)!} \operatorname{H}(x + i c).
\end{equation}

Evaluating our roots at $t=0$ gives $\alpha_1(0) = 1$ and $\alpha_2(0) = 0$. Using $c=1$ for $a_n(x)$ and $c=-1$ for $b_n(x)$, we substitute these formal expansions to obtain:

\begin{equation}
  \begin{aligned}
    a_n(x) & = n [t^n] \mathcal{L}^{-1}\left\{\left(-\alpha_2 s - \ln\left(1 - \frac{\alpha_2 + t}{\alpha_1 + t} e^{zs}\right) - \ln\left(\frac{\alpha_1 + t}{z}\right)\right) s^{-n}\right\}                                              \\
           & = n [t^n] \left( \sum_{i = 1}^{\infty} \frac{(tC)^i}{i} \frac{(x + i z)^{n - 1}}{(n - 1)!} \operatorname{H}(x + i) - \alpha_2 \frac{x^{n-2}}{(n-2)!} - \ln\left(\frac{\alpha_1 + t}{z}\right) \frac{x^{n-1}}{(n-1)!} \right),
  \end{aligned}
\end{equation}

\begin{equation}
  \begin{aligned}
     & b_n(x)                                                                                                                                                                                                                                             \\
     & = n (-1)^n [t^n] \mathcal{L}^{-1}\left\{\left(\alpha_2 s - \ln\left(1 - \frac{\alpha_2 (\alpha_2 + t)}{\alpha_1 (\alpha_1 + t)} e^{-zs}\right) - \ln\left(\frac{\alpha_1 (\alpha_1 + t)}{z}\right)\right) s^{-n} \right\}                          \\
     & = n (-1)^n [t^n] \left( \sum_{i = 1}^{\infty} \frac{(tC)^{3i}}{i} \frac{(x - i z)^{n - 1}}{(n - 1)!} \operatorname{H}(x - i) + \alpha_2 \frac{x^{n-2}}{(n-2)!} - \ln\left(\frac{\alpha_1 (\alpha_1 + t)}{z}\right) \frac{x^{n-1}}{(n-1)!} \right).
  \end{aligned}
\end{equation}

For $c_n(x)$, expanding the rational function yields:

\begin{equation}
  \frac{Q(t, s)}{R(t, s)} = \frac{e^s}{t^3} \left(z \sum_{i=0}^{\infty} (tC)^{3i} e^{-zis} - \alpha_1^2\right),
\end{equation}

\begin{equation}
  \begin{aligned}
    c_n(x) & = (-1)^{n-1} [t^{n-1}] \mathcal{L}^{-1}\left\{\frac{Q(t, s)}{R(t, s)} e^{2s} s^{-n-1}\right\}                                                                                          \\
           & = (-1)^{n-1} \Bigg( 2 [t^{n+1}] \sum_{i=0}^{\infty} (tC)^{3i-1} \frac{(x + 3 - i z)^{n}}{n!} \operatorname{H}(x + 3 - i)                                                               \\
           & \quad - [t^{n+2}] \sum_{i=0}^{\infty} (tC)^{3i} \frac{(x + 3 - i z)^{n}}{n!} \operatorname{H}(x + 3 - i) - [t^{n+2}] \alpha_1^2 \frac{(x + 3)^{n}}{n!} \operatorname{H}(x + 3) \Bigg).
  \end{aligned}
\end{equation}

The core parts of these expressions take the form of $ [t^n] (tC)^p (a + b z)^q $. To extract explicit coefficients, we apply the Lagrange Inversion Theorem via the substitution $u = t^2$. The shifted Catalan function $f(u) = C(u) - 1$ has the compositional inverse $g(u) = \frac{u}{(u+1)^2}$. By shifting to the $u$-domain:

\begin{equation}
  [t^n] (tC)^p (a + b z)^q = [u^{\frac{n-p}{2}}] (u+1)^{n-q-1} \big((a-b)u+(a+b)\big)^q (1-u).
\end{equation}

Surprisingly, we observe that across all of our previous summations for $a_n(x)$, $b_n(x)$, and $c_n(x)$, the exponent $n - q - 1$ consistently evaluates to either $0$ or $1$.

This trivializes the extraction. Because $(u+1)^0 = 1$ and $(u+1)^1 = 1 + u$, the term $[u^{\frac{n-p}{2}}]$ merely extracts the corresponding coefficients from the remaining low-degree polynomial $( (a-b)u + (a+b) )^q (1-u)$. As a result, the seemingly infinite fractional generating functions cleanly collapse into a finite polynomial expansion of basic binomials. We omit the details of extracting the coefficients for the remaining terms.

\section{Exact Closed-Form Probability Distributions}

After applying the Lagrange inversion and simplifying, the raw, dimensionally-consistent Lebesgue measures are given by the explicit piecewise polynomials:

\begin{equation}
  \begin{aligned}
    a_n(x) & = \frac{n}{(n-1)!} \sum_{\substack{i=1                                                                      \\ i \equiv n \pmod 2}}^{n} \frac{1}{i} \operatorname{H}(x + i) \Bigg[ \binom{n-1}{\frac{n-i}{2}} (x-i)^{\frac{n-i}{2}} (x+i)^{\frac{n+i-2}{2}} \\
           & \quad - \binom{n-1}{\frac{n-i-2}{2}} (x-i)^{\frac{n-i-2}{2}} (x+i)^{\frac{n+i}{2}} \Bigg]                   \\
           & \quad - \frac{2^{n-1}}{(n-1)!} x^{n-1} + \frac{1}{2} \binom{n}{\frac{n}{2}} \frac{x^{n-2} (x - n)}{(n-1)!},
  \end{aligned}
\end{equation}

\begin{equation}
  \begin{aligned}
    b_n(x) & = \frac{n (-1)^n}{(n-1)!} \sum_{\substack{i=1                                                                           \\ i \equiv n \pmod 2}}^{\lfloor n/3 \rfloor} \frac{1}{i} \operatorname{H}(x - i) \Bigg[ \binom{n-1}{\frac{n-3i}{2}} (x+i)^{\frac{n-3i}{2}} (x-i)^{\frac{n+3i-2}{2}} \\
           & \quad - \binom{n-1}{\frac{n-3i-2}{2}} (x+i)^{\frac{n-3i-2}{2}} (x-i)^{\frac{n+3i}{2}} \Bigg]                            \\
           & \quad + \frac{(-1)^{n-1} 2^{n-1}}{(n-1)!} x^{n-1} + \frac{1}{2} \binom{n}{\frac{n}{2}} \frac{x^{n-2} (3x + n)}{(n-1)!},
  \end{aligned}
\end{equation}

\begin{equation}
  \begin{aligned}
    c_n(x) & = \frac{(-1)^{n-1}}{n!} \sum_{\substack{i=0                                                             \\ i \equiv n \pmod 2}}^{\lfloor (n+2)/3 \rfloor} \operatorname{H}(x + 3 - i) \Bigg[ \binom{n}{\frac{n-3i+2}{2}} (x+3+i)^{\frac{n-3i+2}{2}} (x+3-i)^{\frac{n+3i-2}{2}} \\
           & \quad - 2\binom{n}{\frac{n-3i}{2}} (x+3+i)^{\frac{n-3i}{2}} (x+3-i)^{\frac{n+3i}{2}}                    \\
           & \quad + \binom{n}{\frac{n-3i-2}{2}} (x+3+i)^{\frac{n-3i-2}{2}} (x+3-i)^{\frac{n+3i+2}{2}} \Bigg]        \\
           & \quad + \frac{2 (-1)^{n-1}}{n+2} \binom{n}{\frac{n}{2}} \frac{(x + 3)^{n}}{n!} \operatorname{H}(x + 3).
  \end{aligned}
\end{equation}

To translate these bounded measures back into geometric probabilities, we divide them by the total unrestricted volume of the uniform probability simplex. For $n$ sequence variables under the constraint $\sum (x_i - 1) = x$, the sum $\sum x_i = x + n$ bounds a total unrestricted volume of $\frac{(x+n)^{n-1}}{(n-1)!}$. Similarly, for $n+1$ variables, the total volume is $\frac{(x+n+1)^n}{n!}$.

Dividing the measures by these respective volumes and mapping the variables via $x = \frac{2}{1 - w} - N$ (for $P_c(N-1; N, w)$), $x = \frac{2}{w} - N$ (for $P_c(3; N, w)$) and $x = \frac{2}{w} - N - 1$ (for $P(3; N, w)$) mathematically resolves the dimensional polynomial factors and normalizes the constraints.

By analyzing the strict geometric bounds of these probabilities, we recognize that the remaining Heaviside step functions serve only to enforce absolute threshold limits. For $P_c(N-1)$, substituting the bounds transforms the expression directly into a finite combinatorial sum bounded by $\lfloor \frac{1}{1 - w} \rfloor$. For $b_n$ and $c_n$, the probabilities inherently equate to zero outside their respective absolute thresholds ($w \ge \frac{2}{N}$ and $w \ge \frac{2}{N-2}$). Within their strictly valid domains, replacing the index $i$ with $N - 2p$ (for $P_c(N-1)$) and $2p - N$ (for $P_c(3)$ and $P(3)$) perfectly collapses the Heaviside inequalities, naturally establishing exactly $\lfloor \frac{1}{w} \rfloor$ (or $\lfloor \frac{1}{w} \rfloor + 1$) as the explicit upper limit for the summations.

Adopting the standard combinatorial extension where $\binom{n}{m} = 0$ if $m < 0$, $n < m$, or $m \notin \mathbb{Z}$, we obtain our final, elegantly reduced exact closed-form formulae:

For $w < 1 - \frac{2}{N}$ (with $P_c(N-1; N, w) = 1$ otherwise), we obtain:
\begin{equation}
  \begin{aligned}
     & 1 - P_c(N-1; N, w)                                                                                                                        \\
     & = \sum_{p = 0}^{\lfloor \frac{1}{1 - w} \rfloor} \binom{N}{p} \big(1-N(1-w)\big) \big(1-p(1-w)\big)^{N-p-1} \big(1-(N-p)(1-w)\big)^{p-1}.
  \end{aligned}
\end{equation}

For $w < \frac{2}{N}$ (with $P_c(3; N, w) = 1$ otherwise), we obtain:
\begin{equation}
  \begin{aligned}
    1 - P_c(3; N, w) & = (-1)^{N-1} \Bigg[ (2 - Nw)^{N - 1} + \frac{1}{2} \binom{N}{\frac{N}{2}} (Nw - 3) \left(1 - \frac{N}{2}w\right)^{N - 2}                                                   \\
                     & \quad + \sum_{p = \lceil \frac{N+1}{2} \rceil}^{\lfloor \frac{1}{w} \rfloor} (Nw - 3) \binom{N}{3p - N} (1 - pw)^{3p - N - 1} \big(1 - (N - p)w\big)^{2N - 3p - 1} \Bigg].
  \end{aligned}
\end{equation}

For $w < \frac{2}{N-2}$ (with $P(3; N, w) = 1$ otherwise), we obtain:
\begin{equation}
  \begin{aligned}
     & 1 - P(3; N, w)                                                                                                                                                                         \\
     & = (-1)^{N-1} \Bigg[ -2 \binom{N}{\frac{N}{2}} \frac{1}{N + 2} \left(1 - \left(\frac{N}{2} - 1\right)w\right)^N                                                                         \\
     & \quad + \sum_{p = \lceil \frac{N+1}{2} \rceil}^{\lfloor \frac{1}{w} \rfloor + 1} \Bigg( \binom{N}{3p - N - 1} \big(1 - (p - 1)w\big)^{3p - N - 1} \big(1 - (N - p - 1)w\big)^{2N-3p+1} \\
     & \qquad - 2\binom{N}{3p - N} \big(1 - (p - 1)w\big)^{3p - N} \big(1 - (N - p - 1)w\big)^{2N-3p}                                                                                         \\
     & \qquad + \binom{N}{3p - N + 1} \big(1 - (p - 1)w\big)^{3p - N + 1} \big(1 - (N - p - 1)w\big)^{2N-3p-1} \Bigg) \Bigg].
  \end{aligned}
\end{equation}

\bibliographystyle{amsplain}
\bibliography{references}

\end{document}